\documentclass{amsart}

\usepackage{citesort}

\newtheorem{theorem}{Theorem}[section]
\newtheorem{lemma}[theorem]{Lemma}

\theoremstyle{definition}

\theoremstyle{remark}
\newtheorem{remark}[theorem]{Remark}

\numberwithin{equation}{section}

\newcommand{\abs}[1]{\lvert#1\rvert}
\def\R{{\mathbb{R}}}

\def\epsilon{{\varepsilon}}
\def\phi{{\varphi}}
\def\theta{{\vartheta}}
\def\dt{{\frac{d}{dt}}}

\DeclareMathOperator{\divergence}{div}
\DeclareMathOperator{\divergenz}{div}

\def\xx{{\frac{x}{\abs x}}}
\def\txx{{\tfrac{x}{\abs x}}}

\begin{document}

\title{Decay at Infinity for Parabolic Equations}

\author{Oliver C. Schn\"urer}
\address{Oliver Schn\"urer, Freie Universit\"at Berlin,
  Ar\-nim\-al\-lee 2--6, 14195 Berlin, Germany}
\def\fudomain{@math.fu-berlin.de}
\email{Oliver.Schnuerer\fudomain}

\author{Hartmut R. Schwetlick}
\address{Mathematical Sciences, University of Bath, Bath BA2 7AY, UK}
\email{H.Schwetlick@bath.ac.uk}


\date{August 2005.}


\begin{abstract}
We consider solutions to linear parabolic equations
with initial data decaying at spatial infinity.
For a class of advection-diffusion equations with a 
spatially dependent velocity field,
we study the behavior of solutions as time tends to infinity.
We characterize velocity fields, so that positive solutions
decay or lift-off at spatial infinity
as time tends to infinity. This addresses
the question of stability of the zero solution for decaying
perturbations.
\end{abstract}

\maketitle

\section{Introduction}
Consider solutions $u:\R^n\times\R_+\to\R$ of the parabolic equation
\begin{equation}\label{flow eqn}
\dot u=\Delta u+\left\langle b,\nabla u\right\rangle\quad
\text{in~}\R^n\times\R_+,
\end{equation}
where $b:\R^n\to\R^n$ is of the form $b(x)=\frac x{\abs{x}}\psi(\abs{x})$,
such that the flow equation preserves rotational symmetry of solutions.
We are interested in the long-time behavior of solutions $u$ that
decay initially, $u(x,0)\to0$ for $\abs x\to\infty$.
For bounded vector fields $b$, a solution lifts off, i.\,e.\
$$\lim\limits_{|x|\to\infty}\lim\limits_{t\to\infty}u(x,\,t)>0,$$
or converges to zero uniformly as $t\to\infty$.
Convergence to zero corresponds to dynamical stability of
the zero solution.
In \cite{OSAlbert}, stability of symmetric gradient K\"ahler-Ricci
solitons was proved for perturbations that decay at infinity.
The equation arising there behaves similarly to \eqref{flow eqn},
when $\left\langle b,x\right\rangle\ge 0$.
We prove in Theorem \ref{abfall thm}
that in this case
bounded solutions $u$ of \eqref{flow eqn}
which decay initially,
must tend to zero as $t\to\infty$.

The situation is different when $\left\langle b,x\right\rangle\le0$. For
$b=-x$, no positive solution tends to zero as $t\to\infty$.
However, it will be shown that linear growth at infinity
is not necessary for the lift-off phenomenon.
We characterize precisely the critical growth rate
for large $x$ for the vector fields
$b(x)=-\frac x{\abs{x}}\psi(\abs{x})$ to be $\psi(r)\approx 1/r$.
In Section \ref{lift-off} we will show that vector fields of faster growth
lead to the lift-off of positive solutions.
On the contrary, we show in Section \ref{conv to zero}
that slower growth forces solutions to converge uniformly  to zero.

The dependence of the behavior of solutions on the sign of
$\left\langle b,x\right\rangle$ can be understood as follows. Rotationally
symmetric solutions that decay monotonically in $\abs{x}$ may serve
as barriers. For these
functions, $\left\langle x,\nabla u\right\rangle\le0$, so
$\left\langle b,\nabla u\right\rangle\le0$
for $\left\langle b,x\right\rangle\ge0$
and $u$ tends faster to zero
than for the heat equation, where $u$
is known to tend to zero.
If $\left\langle b,x\right\rangle\le0$
comparison with the heat equation is not applicable anymore, and
we may expect that large values of
$\left\langle b,\nabla u\right\rangle$ prevent the solution $u$
to decay to zero.

We wish to thank Albert Chau, Klaus Ecker, J\"urgen Jost,
Stefan M\"uller, the Max Planck Institute for Mathematics
in the Sciences, and Free University Berlin for discussions
and support.

\section{Lift-Off}\label{lift-off}

\subsection{Convergence to a Constant}
Throughout that paper we will use the following lemma.
Here, the vector field
$b(x)=\xx\psi(|x|)$ is not assumed
to point in a specific direction.
We remark that we use the space $C^\alpha\left(\R^n\right)$ for
uniformly bounded functions
having bounded H\"older semi-norm with exponent $\alpha$.

\begin{lemma}\label{gen b lem}
Let $\psi:\R_+\to\R$, $0<\alpha<1$ be such that
$b(x)=\frac x{|x|}\psi(|x|)\in C^\alpha\left(\R^n\right)$.
Let $u_0\in C^{2+\alpha}\left(\R^n\right)$. Then there exists a
unique positive solution $u\in C^{2+\alpha,\,1+\alpha/2}
\left(\R^n\times\R_+\right)$ of
\begin{equation}\label{model case}
\begin{cases}
\dot u=\Delta u+\left\langle b,\nabla u\right\rangle &
\text{in~}\R^n\times\R_+,\\
u(\cdot,0)=u_0 & \text{in~}\R^n.
\end{cases}
\end{equation}
If $u_0$ is rotationally symmetric and monotonically
decreasing (or increasing) in radial
direction (i.\,e.\ $[0,\,\infty)\ni\lambda\mapsto u_0(\lambda x)$
is a monotonically decreasing (or increasing) function for all
$x\in\R^n$), then $u(\cdot,\,t)$
shares these properties for any $t>0$.
Moreover, the assumptions above guarantee that $u(\cdot,\,t)$
tends locally uniformly to a constant as $t\to\infty$.
If $u_0\ge0$ but $u_0$ is not identically zero, then
$u(x,\,t)>0$ for any $x\in\R^n$, $t>0$.
\end{lemma}
\begin{proof}
First note that \cite[Theorem 9.2.3]{KrylovBook96}
provides the claimed regularity for all time
and \cite[Theorem 8.11.1]{KrylovBook96} gives uniform
bounds in $C^{2+\alpha,\,1+\alpha/2}(\R^n\times\R_+)$.

It is only here that we use the boundedness of $\psi$. It
seems possible to weaken this hypothesis.
However, already for
bounded functions $\psi$,
we observe both lift-off of solutions as well as
convergence to zero.
Therefore we will not pursue this
issue any further.

Let $u_0$ be rotationally symmetric and $R$ be any
orthogonal transformation on $\R^n$. Then $u(Rx,\,t)$
is another solution to our initial value problem . As
$u(x,\,t)-u(Rx,\,t)$ vanishes at $t=0$, this is
preserved during the evolution \cite{EckerHuiskenInvent}.
Thus $u(\cdot,\,t)$
stays rotationally symmetric during the evolution.

Let $u_0$ be monotonically decreasing in radial direction.
As in \cite[Appendix~A]{OSAlbert}, we obtain that
$u(\cdot,t)$ is radially decreasing for any fixed $t>0$.

Assume now that $u_0$ is rotationally symmetric and monotonically
decreasing in radial direction. If $u_0$ is increasing,
it suffices to consider $-u_0$ as \eqref{model case} is
a linear equation. We wish to show that $u$ tends to a
constant as $t\to\infty$. The following argument is similar
to \cite{OSAlbert,OSMM,AltschulerWu}. Observe that $u(\cdot,\,t)$
attains its maximum at $x=0$. According to the strong maximum
principle, $u(0,\,t)$ is strictly decreasing in time or $u$
is a constant. The maximum principle \cite{EckerHuiskenInvent}
implies that $\inf u(\cdot,\,t)$ is non-decreasing in time.

Let $h:=\lim\limits_{t\to\infty}u(0,\,t)$. If
$\lim\limits_{t\to\infty}u(x,\,t)=h$ for every $x\in\R^n$,
our uniform a priori estimates guarantee that $u(\cdot,\,t)$
converges locally uniformly to $h$ as $t\to\infty$. Otherwise,
we find $x_0\in\R^n$ and a sequence $t_k\to\infty$ such that
$u(x_0,\,t_k)\le h-\epsilon$ for some positive $\epsilon$.
Define $u_k(x,\,t):=u(x,\,t+t_k)$. As $u$ is uniformly bounded
in $C^{2+\alpha,\,1+\alpha/2}\left(\R^n\times\R_+\right)$,
we can extract a subsequence of $u_k$ that converges locally
uniformly in $C^{2,\,1}\left(\R^n\times\R\right)$ to a
solution $w$ of \eqref{model case} in $\R^n\times\R$.
We obtain
$$w(x_0,\,0)\le h-\epsilon<h=w(0,\,0),$$
so $w$ is not constant. The function $w$ attains its maximum
at $w(0,\,t)$ for all $t\in\R$. According to the strong
maximum principle, this is impossible. We deduce that
$u(x,\,t)\to h$ as $t\to\infty$, locally uniformly in $x$.

Applying the strong maximum principle once again yields that
a non-negative solution becomes immediately positive.
\end{proof}

\subsection{Example for Lift-Off}
Before we state our theorem concerning solutions lifting off
at infinity for $t\to\infty$, we wish to investigate the
following model case.

Consider the evolution equation
\begin{equation}\label{x grad}
\begin{cases}
\dot u=\Delta u-\left\langle x,\nabla u\right\rangle &
\text{in~}\R^n\times\R_+,\\
u(\cdot,0)=u_0 & \text{in~}\R^n
\end{cases}
\end{equation}
for some $u_0\in C^{2+\alpha}$.
Let $w$ be the $C^{2+\alpha,\,1+\alpha/2}$-solution to
$$\begin{cases}
\dot w=\Delta w&\text{in~}\R^n\times[0,1],\\
w(\cdot,0)=u_0&\text{in~}\R^n
\end{cases}$$
as in Lemma \ref{gen b lem}.
It is easy to check that
$$u(x,t):=w\left(e^{-t}x,1-\tfrac12e^{-2t}\right)$$ solves
\eqref{x grad}. If $u_0$ is positive, we see that $u$ converges
exponentially fast to a positive constant as $t\to\infty$.
In particular,
this
shows  that solutions to \eqref{x grad}, which decay at
spatial infinity initially, do not necessarily decay at spatial
infinity in the limit $t\to\infty$.
More precisely,
$$\lim\limits_{\abs{x}\to\infty}\lim\limits_{t\to\infty}
u(x,t)$$ may be different from zero. Note that
the order of the limits is important.
We remark that a formal calculation, that can be made precise
for $u$ with good decay at spatial infinity, suggests
that $\left\vert\int_{\R^n}u\right\vert$ increases exponentially in time,
$$\dt\int\limits_{\R^n}u=\int\limits_{\R^n}\Delta u-\left\langle x,
\nabla u\right\rangle=\int\limits_{\R^n}\divergence(\nabla u-x u)+nu
=n\int\limits_{\R^n}u.$$

\subsection{Lift-Off Theorem}
The following result shows that unbounded vector fields $b$ are not necessary
to let solutions lift off at spatial infinity as $t\to\infty$.

\begin{theorem}\label{abhebe thm}
Let $\psi:\R_+\to\R$ be such that
$b(x)=-\txx\psi(|x|)\in C^{\alpha}\left(\R^n\right)$
for some $0<\alpha<1$. We assume that
\begin{equation}\label{growth cond}
\liminf\limits_{r\to\infty}\frac1{\log r}\int\limits_0^r\psi(\rho)d\rho>n.
\end{equation}
Let $0\le u_0\in C^{2+\alpha}(\R^n)$ with $u_0\not\equiv0$.
In addition, we assume that  $u_0$ is rotationally symmetric
and monotonically decreasing in radial direction.
Then the unique positive solution
$u\in C^{2+\alpha,\,1+\alpha/2}(\R^n\times\R_+)$ of
\begin{equation}\label{flow}
\begin{cases}
\dot u=\Delta u-\left\langle\frac{x}{|x|}\psi,\,\nabla u\right\rangle&
\text{in~}\R^n\times\R_+,\\
u(\cdot,\,0)=u_0&\text{in~}\R^n
\end{cases}
\end{equation}
is rotationally symmetric and satisfies
$$\lim\limits_{t\to\infty}\sup\limits_{x\in\R^n}u(x,\,t)=
\inf\limits_{x\in\R^n}\lim\limits_{t\to\infty}u(x,\,t)>0.$$
\end{theorem}
\begin{proof}
Recalling Lemma \ref{gen b lem} it only remains to prove the last convergence claim.
Define $\phi:\R_+\to\R_+$ by
$$\phi(r):=\exp\left(-\int\limits_0^r\psi(\rho)\,d\rho\right).$$
Observe that $\phi$ is bounded and solves $\phi'+\phi\psi=0$.
Define also
$$I_R(t):=\int\limits_{B_R(0)}\phi(|x|)\cdot u(x,\,t)\,dx$$
and
$$I(t):=\int\limits_{\R^n}\phi u.$$
Our assumptions ensure that there exists $\epsilon>0$,
$r_0>0$ such that $\phi(r)\le r^{-n-\epsilon}$
for all $r\ge r_0$.
Thus, we deduce
$$\int\limits_{\R^n}\phi(|x|)dx<\infty,$$
$I_R(0)\le I(0)<\infty$, and $0<I(0)$.
We compute for $\phi u=\phi(|x|)u(x,\,t)$
\begin{align*}
\dt(\phi u)=&\phi\dot u\\
=&\phi\Delta u-\phi\left\langle\txx\psi,\,\nabla u\right\rangle\\
=&\divergenz(\phi\nabla u)-\left\langle\phi'\txx,\,\nabla u\right\rangle
-\phi\left\langle\txx\psi,\,\nabla u\right\rangle\\
=&\divergenz(\phi\nabla u)-(\phi'+\phi\psi)\left\langle\txx,\,
\nabla u\right\rangle\\
=&\divergenz(\phi\nabla u).
\end{align*}
So we obtain that
$$\dt I_R(t)=\int\limits_{\partial B_R}\phi
\left\langle\txx,\,\nabla u\right\rangle\to0\quad\text{for~}R\to\infty,$$
as $\abs{\nabla u}$ is bounded
and $\phi$ decays faster than $r^{-n}$ at infinity.
Thus $I(t)$ is time independent.

The solution $u(\cdot,\,t)$ stays non-negative during the evolution
and tends to a constant as $t\to\infty$, uniformly on compact subsets
of $\R^n$. Since $I(t)$ is time independent and
$\int\limits_{\R^n}\phi(|x|)dx<\infty$, this constant has to be
positive.
\end{proof}

\begin{remark}
The conditions on $\psi$ in Theorem
\ref{abhebe thm} are fulfilled,
if $\psi$ is $C^\alpha$, vanishes in a neighborhood of the origin, and,
for $r\ge r_0>0$, it is of the form $$\psi(r)=A r^{\beta}
\text{\quad with }\begin{cases}A>0,&\beta>-1;\\A>n,&\beta=-1.\end{cases}$$
\end{remark}

\begin{remark}\label{low bar}
If $u_0$ in Theorem \ref{abhebe thm} is not rotationally symmetric,
$u$ also lifts off. At any positive time $\epsilon$,
$u$ is strictly positive. So there exists $\tilde u_0$
fulfilling the assumptions on $u_0$ in
Theorem \ref{abhebe thm} and $\tilde u_0\le u(\cdot,\,\epsilon)$.
Let $\tilde u$ be the solution to \eqref{flow} with
$\tilde u(\cdot,\,0)=\tilde u_0$.
According to the maximum principle, $\tilde u(x,\,t)\le
u(x,\,t+\epsilon)$. As $\tilde u(\cdot,\,t)$ converges
locally uniformly to a positive constant, we obtain that
$$\inf\limits_{x\in\R^n}\liminf_{t\to\infty}u(x,\,t)>0.$$
\end{remark}

\section{Convergence to Zero}
\label{conv to zero}

In the following, we investigate the behavior of solutions
in the sub-critical case, that is, the vector field
$b$ fails to obey
the growth condition \eqref{growth cond}.

\begin{theorem}\label{abfall thm}
Let $\psi:\R_+\to\R$, $0<\alpha<1$ be such that
$b(x)=-\txx\psi(|x|)\in C^{\alpha}\left(\R^n\right)$
and
$$\limsup\limits_{r\to\infty}\frac1{\log r}
\int\limits_0^r\max\,\{\psi(\rho),\,0\}\,d\rho<n.$$
Let $0\le u_0\in C^{2+\alpha}(\R^n)\cap L^1(\R^n)$ with $u_0\not\equiv0$.
In addition, we assume that  $u_0$ is rotationally symmetric
and monotonically decreasing in radial direction.
Then the unique positive solution
$u\in C^{2+\alpha,\,1+\alpha/2}(\R^n\times\R_+)$ of
$$\begin{cases}
\dot u=\Delta u-\left\langle\frac{x}{|x|}\psi,\,\nabla u\right\rangle&
\text{in~}\R^n\times\R_+,\\
u(\cdot,\,0)=u_0&\text{in~}\R^n
\end{cases}$$
is rotationally symmetric and satisfies
$$\lim\limits_{t\to\infty}\sup\limits_{x\in\R^n}u(x,\,t)=0.$$
\end{theorem}
\begin{proof}
According to Lemma \ref{gen b lem}, we only have to show that
$$\lim\limits_{t\to\infty}\sup\limits_{x\in\R^n}u(x,\,t)=0.$$
Rotational symmetry and monotonicity in radial direction imply that
\begin{equation}\label{mon pres}
\left\langle x,\,\nabla u\right\rangle\le0.
\end{equation}

For $\phi:\R_+\to\R_+$ given by
$$\phi(r)=\exp\left(-\int\limits_0^r\psi_+(\rho)\,d\rho\right),$$
we observe that $\phi$ is bounded and solves $\phi'+\phi\psi_+=0$.
Here we used the decomposition
of $\psi$ in its positive and negative part, $\psi=\psi_+-\psi_-$.

Define $I_R(t)$ and $I(t)$ as above in the proof of
Theorem \ref{abhebe thm}.
Our assumptions on $\psi$ and $u_0$ ensure that
$$I_R(0)\le I(0)<\infty.$$
We compute for $\phi u=\phi(|x|)u(x,\,t)$ as above
$$\dt(\phi u)=\divergenz(\phi\nabla u)
-(\phi'+\phi\psi_+)\left\langle\txx,\,
\nabla u\right\rangle+\phi\psi_-\left\langle\txx,\,
\nabla u\right\rangle.$$
Recalling $\phi'+\phi\psi_+=0$ and \eqref{mon pres}, we deduce that
$$\dt(\phi u)\le\divergenz(\phi\nabla u).$$
Using \eqref{mon pres} again, we get
$$\dt I_R(t)\le\int\limits_{\partial B_R}\phi
\left\langle\txx,\,\nabla u\right\rangle\le0.$$
So we obtain for $0\le t_1\le t_2$ the inequality
$$I(t_1)\ge I(t_2).$$
By assumption, there exists $r_0>0$ such that $\phi(r)\ge r^{-n}$
for all $r\ge r_0$.
Thus, we have
$$\int\limits_{\R^n}\phi(|x|)dx=\infty.$$
The solution $u$ stays non-negative during the evolution.
According to Lemma \ref{gen b lem}, the function $u(\cdot,\,t)$
tends to a constant as $t\to\infty$, uniformly on compact subsets
of $\R^n$. As $I(t)$ is non-increasing in time, this constant has to be
zero.
\end{proof}

\begin{remark}
The conditions on $\psi$ in Theorem \ref{abfall thm} are fulfilled,
if $\psi$ is smooth, vanishes in a neighborhood of the origin, and,
for $r\ge r_0>0$, it is of the form
 $$\psi(r)=A r^{\beta}
\text{\quad with }\begin{cases}A\in\R,&\beta<-1;\\A<n,&\beta=-1.\end{cases}$$
It is easy to check directly that the proof of Theorem
\ref{abfall thm} remains valid for a vector field with
$\psi(r)=nr^{-1}$ outside a compact set.
\end{remark}

\begin{remark}
We want to note that our theorems provide a sharp characterization of
the leading order of the growth rate.
The functions $$\psi(r)=\frac1r\left(n+\frac\alpha{\log r}\right)$$
have all critical growth as
$$\liminf\limits_{r\to\infty}\frac1{\log r}
\int\limits_2^r\psi(\rho)d\rho=n.$$
However, we can show that they
lead to lift-off only if $\alpha>1$, whereas $\alpha\le1$
yields decay to zero. This follows from the respective proofs of
the above theorems and a more detailed investigation of the
integrability of $\phi(|x|)$.
\end{remark}

\begin{remark}
Similar to Remark \ref{low bar}, there is also a version of
Theorem \ref{abfall thm} for $u_0$ not being
rotationally symmetric. Here, we may allow $u_0$ to change sign
too. We can
find a barrier $\beta$, such that $\beta\ge u_0\ge-\beta$ and $\beta$
fulfills the conditions on $u_0$ in Theorem \ref{abfall thm}.
As the solution starting with initial datum $\beta$ tends to zero,
the maximum principle implies that $u(\cdot,\,t)$ converges
uniformly to zero as $t\to\infty$.
\end{remark}

\bibliographystyle{amsplain}

\end{document}